\documentclass[12pt,draft,leqno]{article}
\usepackage{amssymb, eucal, latexsym}
\newtheorem{theorem}{Theorem}[section]

\newtheorem{exa}[theorem]{Example}

\newtheorem{pro}[theorem]{Proposition}
\newtheorem{defi}[theorem]{Definition}

\newtheorem{nota}[theorem]{Notation}

\newtheorem{fact}[theorem]{Fact}

\newtheorem{defpro}[theorem]{Definition and Proposition}

\def\p{\varphi}

\def\d{\delta}

\def\k{\kappa}

\def\LAM{\Lambda}

\def\t{\theta}
\def\TE{\Theta}

\def\pl{\varphi_\Lambda}

\def\lra{\longrightarrow}

\def\lr{\leftrightarrow}

\def\sbe{\subseteq}

\def\stm{\setminus}

\def\fa{\forall}

\def\we{\wedge}

\def\bv{\bigvee}

\def\ap{^\prime}
\def\inv{^{-1}}
\def\st{\ |\ }

\def\lll{\ll_l}
\def\llx{\ll_{\rho}}

\def\llxa{\ll_{\rho\ap}}

\def\nll{\not\ll}

\def\card #1{\vert #1 \vert}

\def\CC{{\cal C}}

\def\HLC{{\bf HLC}}
\def\DHLC{{\bf DHLC}}

\def\IPAC{{\bf MVDOpPerLCCon}}

\def\OPLCC{{\bf OpPerLCCon}}

\def\OPALC{{\bf DOpPerLCCon}}

\def\SLCC{{\bf SkePerLCCon}}

\def\SALC{{\bf DSkePerLCCon}}
\def\ISKLC{{\bf InSkeLC}}
\def\SSKLC{{\bf SuSkeLC}}
\def\ISKAL{{\bf DSuSkeLC}}

\def\SSKAL{{\bf DInSkeLC}}
\def\IOPLC{{\bf InOpPerLC}}
\def\SOPLC{{\bf SuOpPerLC}}
\def\IOPAL{{\bf DSuOpPerLC}}
\def\SOPAL{{\bf DInOpPerLC}}

\def\SKLC{{\bf SkeLC}}
\def\SKAL{{\bf DSkeLC}}

\def\SLC{{\bf SkePerLC}}
\def\SAL{{\bf DSkePerLC}}
\def\OLC{{\bf OpLC}}
\def\OAL{{\bf DOpLC}}
\def\OPLC{{\bf OpPerLC}}
\def\OPAL{{\bf DOpPerLC}}

\def\KCon{{\bf KCon}}

\def\IPA{{\bf MVDOpPerLC}}
\def\IA{{\bf MVDSkePerLC}}
\def\IKA{{\bf MVDSkeLC}}
\def\IIA{{\bf MVDSuSkePerLC}}
\def\SIA{{\bf MVDInSkePerLC}}
\def\IAC{{\bf MVDSkePerLCCon}}

\def\IOLC{{\bf InOpLC}}

\def\IOAL{{\bf DSuOpLC}}

\def\SOLC{{\bf SuOpLC}}

\def\SOAL{{\bf DInOpLC}}
\def\IIPA{{\bf MVDSuOpPerLC}}
\def\SIPA{{\bf MVDInOpPerLC}}

\def\ISLC{{\bf InSkePerLC}}

\def\ISAL{{\bf DSuSkePerLC}}

\def\SSLC{{\bf SuSkePerLC}}

\def\SSAL{{\bf DInSkePerLC}}
\def\IOA{{\bf MVDOpLC}}

\def\IISKAL{{\bf MVDSuSkeLC}}
\def\ISSKAL{{\bf MVDInSkeLC}}
\def\IIOAL{{\bf MVDSuOpLC}}
\def\ISOAL{{\bf MVDInOpLC}}

\def\K{{\bf K}}

\def\1{{\bf 1}}
\def\2{\mbox{{\bf 2}}}
\def\3{\mbox{{\bf 3}}}

\def\int{\mbox{{\rm int}}}

\def\cl{\mbox{{\rm cl}}}

\def\doc{\hspace{-1cm}{\em Proof.}~~}
\def\sq{\hspace*{\fill} \hbox{\vrule\vbox{\hrule\phantom{o}\hrule}\vrule}}
\def\sqs{\sq \vspace{2mm}}

\def\BBBB{\mathbb{B}}

\def\bU0{\bar{U}=(U^0,(U^i,U^{ci})_{i\in\omega})}

\def\bV0{\bar{V}=(V^0,(V^i,V^{ci})_{i\in\omega})}

\def\Bool{{\bf BoolAlg}}

\def\doc{\hspace{-1cm}{\em Proof.}~~}
\def\sq{\hspace*{\fill} \hbox{\vrule\vbox{\hrule\phantom{o}\hrule}\vrule}}
\def\sqs{\sq \vspace{2mm}}

\title{{\LARGE\bf
Some Isomorphism Theorems for MVD-algebras}\\
\vspace{0.35cm}
{\large\bf Elza Ivanova-Dimova}\thanks{This paper was supported by the project No. 140/08.05.2014 ,,Function spaces and dualities” of the Sofia University ,,St. Kl. Ohridski”.}\\
\vspace{0.25cm}
 {\footnotesize Dept. of Math. and
Informatics, Sofia University,  J. Bourchier 5, 1126 Sofia,
Bulgaria}}

\date{}

\begin{document}
\maketitle
\begin{abstract}
{\footnotesize
\noindent In the papers \cite{D-TA09,D-AMH1-10,D-AMH2-10} many Stone-type duality theorems for the category of locally compact Hausdorff spaces and continuous maps and some of its subcategories were proved. The dual objects in all these theorems are the local contact algebras. In \cite{VDDB} the notion of an MVD-algebra was introduced and it was shown that it is equivalent to the notion of a local contact algebra. In this paper we express the duality theorems mentioned above in a new form using MVD-algebras and appropriate morphisms between them instead of local contact algebras and the respective morphisms.}
\end{abstract}

{\footnotesize {\em  MSC:} primary 54D45, 18A40; secondary 06E15,
54C10, 54E05.

}

\footnotetext[1]{{\footnotesize {\em E-mail address:}
elza@fmi.uni-sofia.bg}}

\baselineskip = \normalbaselineskip

\section*{Introduction}

The idea of building a {\it region-based theory of space} belongs to A. N. Whitehead \cite{Whitehead} and T. de Leguna \cite{de-Laguna}. Survey papers describing various aspects and historical remarks on region-based theory of space are \cite{Gerla,BeDu,Vak,Pratt}. With the help of the notion of a {\it region-based topology} (which is called {\it local contact algebra} (briefly, $LC$-{\it algebra}) in \cite{DV1}) Roeper \cite{Roeper} gave one of the possible first-order axiomatizations for region-based theory of space. The notion of the region-based topology is based on two primitive spatial relations: {\it contact} and the one-place relation of {\em limitedness}. An attempt to give a different formulation of the same theory using only  one primitive relation, called {\em interior parthood}, was made by Mormann \cite{Mormann} but, as it was pointed out in \cite{VDDB}, the obtained notion  of an {\em enriched Boolean algebra}\/ was more general than it was necessary; the right axiomatization with only one primitive relation was given in \cite{VDDB}, where the appropriate notion of an {\em MVD-algebra}\/ was introduced. In \cite{D-TA09,D-AMH1-10} Dimov defined categories {\bf DHLC}, $\SKAL$, $\SAL$, $\OAL$ and $\OPAL$ whose objects are all complete local contact algebras and whose morphisms
are some appropriate functions between them. These categories are dual to the categories of all locally compact Hausdorff spaces and, respectively, all continuous maps, all continuous skeletal maps, all skeletal perfect maps, all open maps and all open perfect maps. Here we define five categories $\IKA$, $\IA$, $\IOA$, $\IPA$ and {\bf MVDHLC} whose objects are all complete MVD-algebras and whose morphisms are some appropriate functions between them, and we prove that these categories are isomorphic, respectively, to the categories $\SKAL$, $\SAL$, $\OAL$, $\OPAL$ and {\bf DHLC}.

We now fix the notations.

If $\CC$ denotes a category, we write $X\in \card\CC$ if $X$ is
 an object of $\CC$, and $f\in \CC(X,Y)$ if $f$ is a morphism of
 $\CC$ with domain $X$ and codomain $Y$.

All lattices will be with top (= unit) and bottom (= zero)
elements, denoted respectively by 1 and 0. We do not require the
elements $0$ and $1$ to be distinct.



If $(X,\tau)$ is a topological space and $M$ is a subset of $X$,
we denote by $\cl_{(X,\tau)}(M)$ (or simply by $\cl(M)$) the
closure of $M$ in $(X,\tau)$ and by $\int_{(X,\tau)}(M)$ (or
briefly by $\int(M))$ the interior of $M$ in $(X,\tau)$.

The  closed maps and the open maps between topological spaces are
assumed to be continuous but are not assumed to be onto. Recall
that a map is {\em perfect}\/ if it is closed and compact (i.e.
point inverses are compact sets). A continuous map $f:X\lra Y$ is called {\em quasi-open} (\cite{MP}) if for every non-empty open subset $U$ of $X$, $\int(f(U))\not=\emptyset$; a function $f:X\longrightarrow Y$ is called {\em skeletal} 
if $\int(\cl(f(U)))\not=\emptyset$, for every non-empty open subset $U$ of $X$.

\section{Preliminaries}
\begin{defpro}\label{ladj}
\rm Let us recall the notion of {\em lower adjoint} for posets.
Let $\p:A\lra B$ be an order-preserving map between posets. If
$$\pl:B\lra A$$ is an order-preserving map satisfying the following
condition

\medskip

\noindent($\LAM$) for all $a\in A$ and all
$b\in B$, $b\le \p(a)$ iff $\pl (b)\le a$

\medskip

\noindent(i.e.,  the pair
$(\pl,\p)$ forms a {\em Galois connection between posets $B$ and $A$}\/) then we will say that $\pl$ is a {\em lower adjoint}\/ of $\p$.
It is easy to see that condition ($\LAM$) is equivalent to the following  two conditions:

\smallskip

\noindent ($\LAM$1) $\fa b\in B$, $\p(\p_\LAM (b))\ge b$;

\smallskip

\noindent ($\LAM$2) $\fa a\in A$, $\p_\LAM (\p(a))\le a$.

\smallskip

%
%
%
%
%
%






%
%


\end{defpro}

\begin{fact}\label{L2rave}{\rm (\cite{D-TA09})}
 If $A$ and $B$ are Boolean algebras, $\p:A\lra B$ is a
Boolean homomorphism, $A$ has all meets and $\p$ preserves them,
then
%
%
%
%
$\fa a\in A$  and  $\fa b\in B$, $\pl(\p(a)\we
b)=a\we\pl(b)$.
\end{fact}

\begin{defi}\label{conalg}
\rm
An algebraic system $(B,0,1,\vee,\we, {}^*, C)$ is called a {\it
contact Boolean algebra}\/ or, briefly, {\it contact algebra}
(abbreviated as CA or C-algebra) (\cite{DV1})
 if the system
$(B,0,1,\vee,\we, {}^*)$ is a Boolean algebra (where the operation
$``$complement" is denoted by $``\ {}^*\ $")
  and $C$
is a binary relation on $B$, satisfying the following axioms:

\smallskip

\noindent (C1) If $a\not= 0$ then $aCa$;\\
(C2) If $aCb$ then $a\not=0$ and $b\not=0$;\\
(C3) $aCb$ implies $bCa$;\\
(C4) $aC(b\vee c)$ iff $aCb$ or $aCc$.

\smallskip

\noindent We shall simply write $(B,C)$ for a contact algebra. The
relation $C$  is called a {\em  contact relation}. When $B$ is a
complete Boolean algebra, we will say that $(B,C)$ is a {\em
complete contact Boolean algebra}\/ or, briefly, {\em complete
contact algebra} (abbreviated as CCA or CC-algebra). If $D\sbe B$ and $E\sbe B$,
we will write $``DCE$" for $``(\fa d\in D)(\fa e\in E)(dCe)$".

We will say that two C-algebras $(B_1,C_1)$ and $(B_2,C_2)$ are  {\em
CA-isomorphic} iff there exists a Boolean isomorphism $\p:B_1\lra
B_2$ such that, for each $a,b\in B_1$, $aC_1 b$ iff $\p(a)C_2
\p(b)$. Note that in this paper, by a $``$Boolean isomorphism" we
understand an isomorphism in the category $\Bool$ of Boolean algebras and Boolean homomorphisms.

\smallskip

A CA $(B,C)$ is called {\em connected} if it satisfies the following axiom:

\smallskip

\noindent (CON) If $a\not=0,1$ then $aCa^*$.

\smallskip

A contact algebra $(B,C)$ is called a {\it  normal contact Boolean
algebra}\/ or, briefly, {\it  normal contact algebra} (abbreviated
as NCA or NC-algebra) (\cite{deV,F}) if it satisfies the following axioms (we
will write $``-C$" for $``not\ C$"):

\smallskip

\noindent (C5) If $a(-C)b$ then $a(-C)c$ and $b(-C)c^*$ for some $c\in B$;\\
(C6) If $a\not= 1$ then there exists $b\not= 0$ such that
$b(-C)a$.

\smallskip

\noindent A normal CA is called a {\em complete normal contact
Boolean algebra}\/ or, briefly, {\em complete normal contact
algebra} (abbreviated as CNCA or CNC-algebra) if it is a CCA. The notion of
normal contact algebra was introduced by Fedorchuk \cite{F} under
the name {\em Boolean $\d$-algebra}\/ as an equivalent expression
of the notion of compingent Boolean algebra of de Vries (see its definition below). We call
such algebras $``$normal contact algebras" because they form a
subclass of the class of contact algebras and naturally arise in
normal Hausdorff spaces.

For any CA $(B,C)$, we define a binary relation  $``\ll_C $"  on
$B$ (called {\em non-tangential inclusion})  by $``\ a \ll_C b
\leftrightarrow a(-C)b^*\ $". Sometimes we will write simply
$``\ll$" instead of $``\ll_C$". This relation is also known in the literature under the following names: ``well-inside
relation", ``well below'', ``interior parthood'', ``non-tangential
proper part'' or ``deep
inclusion''.
\end{defi}

The relations $C$ and $\ll$ are inter-definable. For example,
normal contact algebras could be equivalently defined (and exactly
in this way they were introduced (under the name of {\em
compingent Boolean algebras}) by de Vries in \cite{deV}) as a pair
of a Boolean algebra $B=(B,0,1,\vee,\we,{}^*)$ and a binary
relation $\ll$ on $B$ subject to the following axioms:

\smallskip

\noindent ($\ll$1) $a\ll b$ implies $a\leq b$;\\
($\ll$2) $0\ll 0$;\\
($\ll$3) $a\leq b\ll c\leq t$ implies $a\ll t$;\\
($\ll$4) $a\ll c$ and $b\ll c$ implies $a\vee b\ll c$;\\
($\ll$5) If  $a\ll c$ then $a\ll b\ll c$  for some $b\in B$;\\
($\ll$6) If $a\neq 0$ then there exists $b\neq 0$ such that $b\ll
a$;\\
($\ll$7) $a\ll b$ implies $b^*\ll a^*$.

\smallskip

 The proof of
the equivalence of the two definitions of normal contact algebras is
straightforward and analogous to the corresponding statement for proximity
spaces (see Theorems 3.9 and 3.11 in \cite{NW}). One has just to show that
$xCy$ iff $x\nll y^*$.

\smallskip

Obviously, contact algebras could be equivalently defined as a
pair of a Boolean algebra $B$ and a binary relation $\ll$ on $B$
subject to the  axioms ($\ll$1)-($\ll$4) and ($\ll$7).

\smallskip

It is easy to see that axiom (C5) (resp., (C6)) can be stated
equivalently in the form of ($\ll$5) (resp., ($\ll$6)).

\begin{defi}\label{locono}{\rm (\cite{Roeper,DV1})}
\rm An algebraic system $\underline {B}_l=(B,0,1,\vee,\we, {}^*,
\rho, \BBBB)$ is called a {\it local contact algebra} (abbreviated
as LCA)   if $(B,0,1, \vee,\we, {}^*)$ is a Boolean algebra,
$\rho$ is a binary relation on $B$ such that $(B,\rho)$ is a CA,
and $\BBBB$ is a subset of $B$,
satisfying the following axioms:
%
%
%
%

\smallskip

\noindent (BB1) $0\in\BBBB$;\\
(BB2) For $a,b\in B$, $a\le b$ and $b\in\BBBB$ implies $a\in\BBBB$;\\
(BB3) $a,b\in \BBBB$ implies $a\vee b\in\BBBB$;\\
(BC1) If $a\in\BBBB$, $c\in B$ and $a\ll_\rho c$ then there exists
$b\in\BBBB$ such that $a\ll_\rho b\ll_\rho c$ (see \ref{conalg}
for
$``\ll_\rho$");\\
(BC2) If $a\rho b$ then there exists an element $c$ of $\BBBB$
such that
$a\rho (c\we b)$;\\
(BC3) If $a\neq 0$ then there exists  $b\in\BBBB\stm\{0\}$ such
that $b\ll_\rho a$.

\smallskip

Usually, we shall write simply $(B, \rho,\BBBB)$ for a local
contact algebra.  We will say that the elements of $\BBBB$ are
{\em bounded} and the elements of $B\stm \BBBB$  are  {\em
unbounded}. When $B$ is a complete Boolean algebra, we will say
that $(B,\rho,\BBBB)$ is a {\em complete local contact algebra}
(abbreviated by CLCA).

We will say that two local contact algebras $(B,\rho,\BBBB)$ and
$(B_1,\rho_1,\BBBB_1)$ are  {\em LCA-isomorphic} iff there exists
a Boolean isomorphism $\p:B\lra B_1$ such that, for $a,b\in B$,
$a\rho b$ iff $\p(a)\rho_1 \p(b)$ and $\p(a)\in\BBBB_1$ iff
$a\in\BBBB$.

An LCA $(B,\rho,\BBBB)$ is called {\em connected}\/ if the CA $(B,\rho)$ is
connected.
\end{defi}


\begin{exa}\label{rct}
\rm Recall that a subset $F$ of a topological space $(X,\tau)$ is
called {\em regular closed}\/ if $F=\cl(\int (F))$. Clearly, $F$
is regular closed iff it is the closure of an open set.

For any topological space $(X,\tau)$, the collection $RC(X,\tau)$
(we will often write simply $RC(X)$) of all regular closed subsets
of $(X,\tau)$ becomes a complete Boolean algebra
$(RC(X,\tau),0,1,\wedge,\vee,{}^*)$ under the following operations:
$ 1 = X,  0 = \emptyset, F^* = \cl(X\setminus F), F\vee G=F\cup G,
F\wedge G =\cl(\int(F\cap G)).
$
The infinite operations are given by the  formulas:
$\bigvee\{F_{\gamma}\st \gamma\in\Gamma\}=\cl(\bigcup\{F_{\gamma}\st
\gamma\in\Gamma\}),$ and
$\bigwedge\{F_{\gamma}\st \gamma\in\Gamma\}=\cl(\int(\bigcap\{F_{\gamma}\st
\gamma\in\Gamma\})).$

It is easy to see that setting $F \rho_{(X,\tau)} G$ iff $F\cap
G\not = \emptyset$, we define a contact relation $\rho_{(X,\tau)}$ on
$RC(X,\tau)$; it is called a {\em standard contact relation}. So,
$(RC(X,\tau),\rho_{(X,\tau)})$ is a CCA (it is called a {\em
standard contact algebra}). We will often write simply $\rho_X$
instead of $\rho_{(X,\tau)}$. Note that, for $F,G\in RC(X)$,
$F\ll_{\rho_X}G$ iff $F\subseteq \int_X(G)$.

Clearly, if $(X,\tau)$ is a normal Hausdorff space then the
standard contact algebra $(RC(X,\tau),\rho_{(X,\tau)})$ is a
complete NCA.

\end{exa}

In \cite{VDDB} the following notion was introduced:

\begin{defi}\label{MVA3}{\rm (\cite{VDDB})}
\rm A triple $(B,\le,\ll)$ is called an {\em MVD-algebra} if
$(B,\le)$ is a Boolean algebra and the axioms ($\ll$1)-($\ll$6)
(see \ref{conalg})
 as well as the following
two axioms are satisfied:

\smallskip

\noindent ($\ll 4^*$)
\ $a\ll b$ and $a\ll c$ imply $a\ll b\we c$, and

\smallskip

\noindent (MVD) If $a\ll 1$ and $b^*\ll a^*$  then $a\ll b$.

\smallskip

When $(B,\le)$ is a complete Boolean algebra, we will say that
$(B,\le,\ll)$ is a {\em complete MVD-algebra}.
\end{defi}

It follows immediately from the corresponding definitions that
normal contact algebras coincide with MVD-algebras satisfying the
additional axiom

\smallskip

\noindent ($\ll 2\ap$) $1\ll 1$.

\begin{pro}\label{standmv3}{\rm (\cite{VDDB})}
Let $L$ be a locally compact Hausdorff space. Then
$$(RC(L),\sbe,\ll_L),$$
where, for all $F,G\in RC(L)$, $F\ll_L G$ iff $F$ is compact and
$F\sbe \int(G)$, is an MVD-algebra. All such MVD-algebras will be
called\/ {\em standard MVD-algebras}.
\end{pro}

\begin{theorem}\label{eqv3}{\rm (\cite{VDDB})}
 The notions of local contact algebra and MVD-algebra are equivalent.
 More precisely: let $\k$ be the correspondence which assigns to
 every LCA $(B,\rho,\BBBB)$ an MVD-algebra
 $\k(B,\rho,\BBBB)=(B,\le_l,\lll)$, where $a\le_l b$ iff
 $a\we b=a$, and
\begin{equation}\label{mvddefi3}
a\lll b \mbox{ iff } a\in\BBBB \mbox{ and } a\llx b
\end{equation}
(see \ref{conalg} for $``\ll_\rho$"); further, let $\t$ be the
correspondence which assigns to
 each  MVD-algebra
 $(B,\le,\ll)$   an LCA $(B,\rho_m,\BBBB_m)=\t(B,\le,\ll)$,
where
\begin{equation}\label{defbb3}
\BBBB_m=\{a\in B\st a\ll 1\}
\end{equation}
and, for $a,b\in B$,
\begin{equation}\label{fx3}
a\ll_{\rho_m} b\lr(\fa c\ll 1)[(c\we a)\ll(c^*\vee b)].
\end{equation}
(or, equivalently, $a\rho_m b$  iff there exists   $c\ll 1$  such
that  $(c\we a)\nll (c\we b)^*)$. Then $\k$ and $\t$ are bijective
correspondences between the classes of all LCA's and all
MVD-algebras, and $\k=\t\inv$.
\end{theorem}

The following obvious fact was noted in
\cite{Biacino-and-Gerla-1996}.

\begin{fact}\label{confact}{\rm (\cite{Biacino-and-Gerla-1996})}
Let $(X,\tau)$ be a topological space. Then the standard contact
algebra $(RC(X,\tau),\rho_{(X,\tau)})$ is connected iff the space
$(X,\tau)$ is connected.
\end{fact}




%




\begin{pro}\label{semsk}
(a) Every quasi-open map is skeletal.

\noindent (b) {\rm  (\cite{D-TA09})} Let $X$ be a regular space and $f:X\lra Y$ be
a closed map. Then $f$ is quasi-open iff $f$ is skeletal.
\end{pro}

\begin{defi}\label{bounded}{\rm  (\cite{D-TA09})}
\rm
Let $(B,\rho,\BBBB)$ be a local contact algebra.
An ultrafilter $u$ in $B$ is called a {\em bounded ultrafilter} if $u\cap\BBBB\not=\emptyset$.
\end{defi}

\begin{nota}\label{In&Su}
\rm If ${\bf K}$ is a category, then by ${\bf InK}$ (resp., ${\bf SuK}$) we will denote the category having the same objects as the category ${\bf K}$ and whose morphisms are only the injective (resp., surjective) morphisms of ${\bf K}$
\end{nota}

\begin{nota}\label{convcatn}
\rm If $\K$ is a category whose objects form a subclass of the
class of all topological spaces (resp., contact algebras) then we
will denote by $\KCon$ the full subcategory of $\K$ whose objects
are all $``$connected" $\K$-objects, where $``$connected" is
understood in the usual sense when the objects of $\K$ are
topological spaces and in the sense of \ref{conalg} (see the
condition (CON) there) when the objects of $\K$ are contact
algebras.




\end{nota}

\section{Isomorphism theorems for MVD-algebras}

In  \cite{D-TA09}, a category $\SKAL$ was introduced, namely, the
objects of the category $\SKAL$ are all complete local contact algebras and its morphisms $\varphi:(A,\rho,\mathbb{B})\longrightarrow (B,\eta,\mathbb{B}')$ are all complete Boolean homomorphisms satisfying the following conditions:

\smallskip

\noindent (L1) $\forall a,b\in A$, $\varphi(a)\eta\varphi(b)$ implies $a\rho b$;

\smallskip

\noindent (L2) $b\in\BBBB\ap$ implies $\pl(b)\in\BBBB$ (see \ref{ladj} for $\pl$).

\smallskip

As it was proved in \cite{D-TA09}, the category $\SKAL$ is dually equivalent to the category   $\SKLC$ of all locally compact Hausdorff spaces and all continuous skeletal maps between them.

Let us note that (L1) is equivalent to the following condition:

\smallskip





\noindent (L1') $\forall a,b\in A$, $a\ll_{\rho}b$ implies $\varphi(a)\ll_{\eta}\varphi(b)$.

\begin{defi}\label{catmvdnk}
\rm
Let us define a category which will be denoted by $\IKA$. Its
objects are all complete MVD-algebras (see \ref{MVA3}). If
$(B,\le,\ll)$ and $(B\ap,\le\ap,\ll\ap)$ are two complete
MVD-algebras then $\p:(B,\le,\ll)\lra (B\ap,\le\ap,\ll\ap)$ will
be an $\IKA$-morphism iff $\p:(B,\le)\lra (B\ap,\le\ap)$ is a
complete Boolean homomorphism  satisfying the following axioms:

\smallskip

\noindent (S1) For every $a,b\in B$, [($\fa c\in B$ with $c\ll 1$),
($c\we a\ll c^*\vee b$)] implies [($\fa d\in B\ap$ with $d\ll\ap 1$),
($d\we \p(a)\ll\ap d^*\vee\p(b)$)];\\
(S2) For all $b\in B\ap$, $b\ll\ap 1$ implies  $\pl(b)\ll 1$.

\smallskip

Let the composition of two $\IKA$-morphisms be the usual composition of functions.

It is easy to see that in such a way we have defined a category.
\end{defi}

\begin{theorem}\label{mvdeqvnk}
The categories $\SKAL$ and $\IKA$ are isomorphic; hence the categories $\SKLC$ and $\IKA$ are dually equivalent.
\end{theorem}

\doc Let us define two covariant functors $K:\SKAL\lra \IKA$ and
$\TE:\IKA\lra \SKAL$.

For every $(B,\rho,\BBBB)\in\card\SKAL$ we put
$K(B,\rho,\BBBB)=\k(B,\rho,\BBBB)$ (see \ref{eqv3} for $\k$). Then
Theorem \ref{eqv3} implies that $K$ is well-defined on the objects
of the category $\SKAL$.

Let $\p\in\SKAL((B,\rho,\BBBB),(B\ap,\rho\ap,\BBBB\ap))$. We  will
prove that the same function $\p:B\lra B\ap$ is an $\IKA$-morphism
between $K(B,\rho,\BBBB)$ and $K(B\ap,\rho\ap,\BBBB\ap)$. Since
$\p$ is a complete Boolean homomorphism between Boolean algebras
$B$ and $B\ap$, we need only to check that $\p$ satisfies axioms
(S1) and (S2). Using \ref{eqv3} and (L1'), this can be easily done. So, we
can define:
$$K(\p)=\p.$$
Then, obviously, $K:\SKAL\lra \IKA$ is a (covariant) functor.

Let $(B,\le,\ll)\in\card\IKA$. We put
$\TE(B,\le,\ll)=\t(B,\le,\ll)$ (see \ref{eqv3} for  $\t$). Then
\ref{eqv3} implies that $\TE$ is well-defined on the objects of
the category $\IKA$.

Let $\p\in\IKA((B,\le,\ll),(B\ap,\le\ap,\ll\ap))$. We  will show
that the same function $\p:B\lra B\ap$ is an $\SKAL$-morphism
between $\TE(B,\le,\ll)$ and $\TE(B\ap,\le\ap,\ll\ap)$. For doing
this it is enough to prove that $\p$ satisfies conditions (L1) and
(L2). Using \ref{eqv3} and (L1'), this can be easily done. So, we can
define:
$$\TE(\p)=\p.$$
Then, obviously, $\TE:\IKA\lra \SKAL$ is a (covariant) functor.

From the definition of the functors $K$ and $\TE$ and the
equalities $\k\circ\t=id$, $\t\circ\k=id$ (see \ref{eqv3}), we
conclude that $K\circ\TE=Id_{\IKA}$ and $\TE\circ K=Id_{\SKAL}$.
Hence, the categories $\SKAL$ and $\IKA$ are isomorphic. \sqs

In \cite{D-TA09} a category $\SAL$ was introduced, namely, the objects of the category $\SAL$ are all complete local
contact algebras (see \ref{locono}) and its morphisms $\p:(A,\rho,\BBBB)\lra (B,\eta,\BBBB\ap)$ are all $\SKAL$-morphisms satisfying the following condition:

\smallskip

\noindent(L3) $a\in\BBBB$ implies $\p(a)\in\BBBB\ap$.

\smallskip

Obviously, $\SAL$ is a
subcategory of the category $\SKAL$.

As it was proved in \cite{D-TA09}, the category $\SAL$ is dually equivalent to the category $\SLC$ of all locally compact Hausdorff spaces and all skeletal perfect maps between them.

Note that, by \ref{semsk}(b), the morphisms of the category $\SLC$ are precisely the  quasi-open  perfect maps  (because the perfect maps are closed maps).

\begin{defi}\label{catmvdn}
\rm
Let's define a category which will be denoted by $\IA$. Its
objects are all complete MVD-algebras (see \ref{MVA3}). If
$(B,\le,\ll)$ and $(B\ap,\le\ap,\ll\ap)$ are two complete
MVD-algebras then $\p:(B,\le,\ll)\lra (B\ap,\le\ap,\ll\ap)$ will
be an $\IA$-morphism iff $\p:(B,\le)\lra (B\ap,\le\ap)$ is a
complete Boolean homomorphism  satisfying the axiom (S2) from
\ref{catmvdnk} and the following two additional axioms:

\smallskip

\noindent (ES1) For every $a,b\in B$,  $a\ll b$ implies
$\p(a)\ll\ap\p(b)$;\\
(S3) For all $a\in B$, $a\ll 1$ implies
 $\p(a)\ll\ap 1$.

\smallskip

Let the composition of two $\IA$-morphisms be the usual composition of functions.

It is easy to see that in such a way we have defined a (non-full)
subcategory of the category $\IKA$.
\end{defi}

\begin{theorem}\label{mvdeqvn}
The categories $\SAL$ and $\IA$ are isomorphic; hence the categories $\SLC$ and $\IA$ are dually equivalent.
\end{theorem}

\doc We will show that the restrictions $K_p:\SAL\lra \IA$ and
$\TE_p:\IA\lra \SAL$ of the functors $K:\SKAL\lra \IKA$ and
$\TE:\IKA\lra\SKAL$ defined in the proof of Theorem \ref{mvdeqvnk}
are the desired isomorphism functors.

Let $\p\in\SAL((B,\rho,\BBBB),(B\ap,\rho\ap,\BBBB\ap))$. Then, as
it was shown in the proof of \ref{mvdeqvnk},  the same function
$\p:B\lra B\ap$ is an $\IKA$-morphism between
MVD-algebras
$K(B,\rho,\BBBB)$ and $K(B\ap,\rho\ap,\BBBB\ap)$. So, we need only
to check that $\p$ satisfies axioms (ES1) and (S3).

Put $K_p(B,\rho,\BBBB)=(B,\le,\ll)$ and
$K_p(B\ap,\rho\ap,\BBBB\ap)=(B\ap,\le\ap,\ll\ap)$. Then, by
\ref{eqv3}, $a\ll b$ iff $a\in\BBBB$ and $a\llx b$; also $a\ll\ap
b$ iff $a\in\BBBB\ap$ and $a\llxa b$. Using (L3), we get easily
that (S3) is fulfilled. We will show that (ES1) takes place.
Let $a,b\in B$ and  $a\ll b$. Then $a\llx b$ and $a\in\BBBB$.
Thus, by (L1'), $\p(a)\llxa\p(b)$. Since, by (L3),
$\p(a)\in\BBBB\ap$, we obtain that $\p(a)\ll\ap\p(b)$. We have
established that
$\p\in\IA(K_p(B,\rho,\BBBB),K_p(B\ap,\rho\ap,\BBBB\ap))$.

Let $\p\in\IA((B,\le,\ll),(B\ap,\le\ap,\ll\ap))$. We will show
that $\p$ satisfies condition (S1). Let $a,b\in B$ and let, for
every $c\in B$ with $c\ll 1$, $c\we a\ll c^*\vee b$ holds. Take
$d\in B\ap$ with $d\ll\ap 1$. Then, by (S2), $c=\pl(d)\ll 1$.
Hence $c\we a\ll c^*\vee b$. Using (ES1), we obtain that
$\p(c)\we\p(a)\ll\ap (\p(c))^*\vee \p(b)$. Then ($\LAM$1) and
($\ll$3) imply that $d\we\p(a)\ll\ap d^*\vee\p(b)$. Hence (S1) is
established. Now, as it was shown in the proof of \ref{mvdeqvnk},
the same function $\p:B\lra B\ap$ is an $\SKAL$-morphism between
$\TE(B,\le,\ll)$ and $\TE(B\ap,\le\ap,\ll\ap)$.  So, we need only
to prove that $\p$ satisfies condition (L3). This can be done
readily, using (S3).
 The rest follows from Theorem \ref{mvdeqvnk}. \sqs

In \cite{D-TA09} a category $\OAL$ was introduced, namely, the objects of the category $\OAL$ are all complete local
contact algebras (see \ref{locono}) and its morphisms $\p:(A,\rho,\BBBB)\lra (B,\eta,\BBBB\ap)$ are all $\SKAL$-morphisms satisfying the following condition:

\smallskip

\noindent (LO) $\fa a\in A$ and $\fa b\in \BBBB\ap$,
$\p_\LAM(b)\rho a$ implies
 $b\eta\p(a)$.

\smallskip

Obviously, $\OAL$ is a (non-full) subcategory of the category $\SKAL$.

As it was proved in \cite{D-TA09}, the category $\OLC$ of all locally compact Hausdorff spaces and all open maps between them and the category $\OAL$ are dually equivalent.

\begin{defi}\label{catmvdo}
\rm Let us define a category which will be denoted by $\IOA$. Its
objects are all complete MVD-algebras. If $(B,\le,\ll)$ and
$(B\ap,\le\ap,\ll\ap)$ are two complete MVD-algebras then
$\p:(B,\le,\ll)\lra (B\ap,\le\ap,\ll\ap)$ will be an
$\IOA$-morphism iff $\p$ is an $\IKA$-morphism (see
\ref{catmvdnk}) which  satisfies the following axiom:

\smallskip

\noindent (SO) For all $a\in B$ and $b\in B\ap$, $b\ll\ap\p(a)$
implies
 $\p_{\LAM}(b)\ll a$
(see \ref{ladj} for  $\pl$).

\smallskip

\noindent Let the composition of two $\IOA$-morphisms be the usual composition of functions.
It is easy to see that in such a way we have defined a
category.
\end{defi}

\begin{theorem}\label{mvdeqvo}
The categories $\OAL$ and $\IOA$ are isomorphic; hence the categories $\OLC$ and $\IOA$ are dually equivalent.
\end{theorem}

\doc  We will show that the restrictions $K_o:\OAL\lra \IOA$ and
$\TE_o:\IOA\lra \OAL$ of the functors $K$ and $\TE$ defined in the
proof of Theorem \ref{mvdeqvnk} are the desired isomorphism
functors.

Let $\p\in\OAL((B,\rho,\BBBB),(B\ap,\rho\ap,\BBBB\ap))$. We  will
prove that the same function $\p:B\lra B\ap$ is an $\IOA$-morphism
between $K_o(B,\rho,\BBBB)$ and $K_o(B\ap,\rho\ap,\BBBB\ap)$.
Since $\p$ is an $\IKA$-morphism (see the proof of Theorem
\ref{mvdeqvnk}), we need only to check that $\p$ satisfies the
axiom (SO).

Put $K(B,\rho,\BBBB)=(B,\le,\ll)$ and
$K(B\ap,\rho\ap,\BBBB\ap)=(B\ap,\le\ap,\ll\ap)$. Then, by
\ref{eqv3}, $a\ll b$ iff $a\in\BBBB$ and $a\llx b$; also $a\ll\ap
b$ iff $a\in\BBBB\ap$ and $a\llxa b$.

For verifying (SO), note first that (LO) can be formulated
equivalently as:
\begin{equation}\label{l2}
\fa a\in A\mbox{ and }\fa b\in\BBBB\ap,\ b\llxa \p(a) \mbox{ implies }\pl(b)\llx a.
\end{equation}
 Let now $a\in B$, $b\in B\ap$ and $b\ll\ap\p(a)$.
Then $b\llxa \p(a)$ and $b\in \BBBB\ap$. Hence, by (LO),
$\pl(b)\llx a$. Since, by (L2), $\pl(b)\in\BBBB$, we obtain that
$\pl(b)\ll a$.

Therefore, the functor $K_o$ is well-defined.

Let $\p\in\IOA((B,\le,\ll),(B\ap,\le\ap,\ll\ap))$. We  will show
that the same function $\p:B\lra B\ap$ is an $\OAL$-morphism
between $\TE_o(B,\le,\ll)$ and $\TE_o(B\ap,\le\ap,\ll\ap)$. Since,
by the proof of Theorem \ref{mvdeqvnk}, $\p$ is an
$\SKAL$-morphism, it is enough to show that $\p$ satisfies
condition (LO).

Put $\TE(B,\le,\ll)=(B,\rho,\BBBB)$ and
$\TE(B\ap,\le\ap,\ll\ap)=(B\ap,\rho\ap,\BBBB\ap)$.

Let $a\in B$, $b\in \BBBB\ap$ and $b\llxa\p(a)$. Then, by Theorem
(\ref{eqv3}), $b\ll\ap\p(a)$. Hence, by (SO), $\pl(b)\ll a$. This
implies that $\pl(b)\llx a$. Hence, $\p$ satisfies condition (LO).
So, $\TE_o$ is well-defined.

The rest follows from Theorem \ref{mvdeqvnk}.
 \sqs

In \cite{D-TA09}, a category $\OPAL$ was introduced, namely, the objects of the category $\OPAL$ are all complete local
contact algebras (see \ref{locono}) and its morphisms $\p:(A,\rho,\BBBB)\lra (B,\eta,\BBBB\ap)$ are all
$\SAL$-morphisms  satisfying  condition (LO).

Clearly, $\OPAL$ is a subcategory of the category $\SAL$.

As it was proved in \cite{D-TA09}, the category $\OPLC$ all locally compact Hausdorff spaces and all open perfect maps between them and the category $\OPAL$ are dually equivalent.

\begin{defi}\label{catmvd}
\rm Let us now define a subcategory $\IPA$ of the category $\IA$. Its
objects are all complete MVD-algebras. If $(B,\le,\ll)$ and
$(B\ap,\le\ap,\ll\ap)$ are two complete MVD-algebras then a $\IA$-morphism
$\p:(B,\le,\ll)\lra (B\ap,\le\ap,\ll\ap)$ (see \ref{catmvdn}) will be a $\IPA$-morphism if
it satisfies the  axiom (SO) (see \ref{catmvdo}).\\
It is easy to see that in such a way we have defined a category.
\end{defi}

\begin{theorem}\label{mvdeqv}
The categories $\OPAL$ and $\IPA$ are isomorphic; hence the categories $\OPLC$ and $\IPA$ are dually equivalent.
\end{theorem}

\doc  It follows from the proofs of Theorem \ref{mvdeqvn} and Theorem \ref{mvdeqvo}.
 \sqs

In \cite{D-AMH2-10},
a category $\SSKAL$ was introduced, namely,
the objects of the category $\SSKAL$ are  all complete local
contact algebras (see \ref{locono}) and, for any two CLCA's
$(A,\rho,\BBBB)$ and $(B,\eta,\BBBB\ap)$, $\p:(A,\rho,\BBBB)\lra
(B,\eta,\BBBB\ap)$ is an $\SSKAL$-{\em morphism}\/ iff $\p$ is an
$\SKAL$-morphism
which satisfies the following
condition:

\smallskip

\noindent (LS) $\fa a,b\in \BBBB\ap$, $\pl(a)\rho\pl(b)$ implies
$a\eta b$ (see \ref{ladj} for $\pl$).

\smallskip

As it was shown in \cite{D-AMH2-10},
the category  $\SSKAL$ is dually equivalent to the category $\ISKLC$.

Note that condition (LS) is equivalent to the condition below:

\smallskip

\noindent (LS') $\fa a,b\in\BBBB\ap$, $a\ll_{\eta}b$ implies $\pl(a)\ll_{\rho}(\pl(b^*))^*$.

\begin{defi}\label{ISSKAL}
\rm Let $\ISSKAL$ be the category having as objects all complete $MVD$-algebras and let for any two complete $MVD$-algebras $(B,\leq,\ll)$ and $(B',\leq',\ll')$, $\varphi:(B,\leq,\ll)\longrightarrow(B',\leq',\ll')$ be an $\ISSKAL$-morphism iff $\varphi$ is an $\IKA$-morphism (see \ref{catmvdnk}) which satisfies the following condition:

\smallskip

\noindent $(LS'')$ $\forall a,b\in B'$ such that $a,b\ll' 1$, $a\ll' b$ implies $\pl(a)\ll(\pl(b^*))^*$.
\end{defi}

\begin{theorem}\label{izesskal}
The categories $\SSKAL$ and $\ISSKAL$ are isomorphic; hence the categories $\ISKLC$ and $\ISSKAL$ are dually equivalent.
\end{theorem}

\doc We will show that the restrictions $K_r:\SSKAL\lra \ISSKAL$ and $\TE_r:\ISSKAL\lra \SSKAL$ of the functors $K:\SKAL\lra \IKA$ and $\TE:\IKA\lra\SKAL$ defined in the proof of Theorem \ref{mvdeqvnk} are the desired isomorphism functors.

Let $\varphi\in\SSKAL((B,\rho,\BBBB),(B',\rho',\BBBB'))$. We will show that $\varphi:B\longrightarrow B'$ is an $\ISSKAL$-morphism between $MVD$-algebras $K(B,\rho,\BBBB)$ and $K(B',\rho',\BBBB')$. We need only to check that $\varphi$ satisfies $(LS'')$. Let $a,b\in B'$, $a,b\ll'1$ and $a\ll' b$. Then by \ref{eqv3}, $a,b\in\BBBB'$ and $a\ll_{\rho} b$. Since $\varphi$ satisfy (LS'), $\pl(a)\ll_{\rho}(\pl(b^*))^*$. It follows from (S2), that $\pl(a)\ll 1$, i.e. $\pl(a)\ll(\pl(b^*))^*$. Hence $\varphi$ is an $\ISSKAL$-morphism.

Let $\varphi\in\ISSKAL((B,\leq,\ll),(B',\leq',\ll))$. We will show that the same function $\varphi:B\longrightarrow B'$ is a $\SSKAL$-morphism between $\TE(B,\leq,\ll)$ and $\TE(B',\leq',\ll)$. Since, by \ref{mvdeqvnk}, we have that $\varphi$ is an $\SKAL$-morphism, we have only to show that $\varphi$ satisfies $(LS')$. Let $a,b\in\BBBB'$ and $a\ll_{\eta}b$. By, \ref{eqv3} $a,b\ll'1$ and $a\ll'b$. Since $\varphi$ satisfies $(LS'')$, we get that $\pl(a)\ll'(\pl(b^*))^*$, i.e. $\pl(a)\in\BBBB$ and $\pl(a)\ll_{\rho} (\pl(b^*))^*$. Then $\varphi$ satisfy (LS'). Thus $\varphi$ is a $\SSKAL$-morphism. \sqs

In \cite{D-AMH2-10}, a category $\ISKAL$ was introduced, namely, the objects of the category $\ISKAL$ are all complete local contact algebras (see \ref{locono}) and its morphisms $\p:(A,\rho,\BBBB)\lra(B,\eta,\BBBB\ap)$, where $(A,\rho,\BBBB)$ and $(B,\eta,\BBBB\ap)$ are CLCA's, are all $\SKAL$-morphism
which satisfy the following condition:

\smallskip

\noindent (IS) For every bounded ultrafilter $u$ in
$(A,\rho,\BBBB)$ there exists a bounded ultrafilter $v$ in
$(B,\eta,\BBBB\ap)$ such that $\pl(v)\rho u$ (see \ref{bounded},
\ref{ladj} and \ref{conalg} for the notations).

\smallskip

As it was proved in \cite{D-AMH2-10}, the categories $\SSKLC$ and $\ISKAL$ are dually equivalent.

\begin{defi}\label{iiskal}
\rm Let's define a category which will be denoted by $\IISKAL$. Its objects are all complete MVD-algebras. If $(B,\leq,\ll)$ and $(B',\leq',\ll')$ are two complete MVD-algebras then $\varphi\in\IISKAL((B,\leq,\ll),(B',\leq',\ll'))$ iff $\varphi$ is an $\IKA$-morphism (see \ref{catmvdnk}) which satisfies the following axiom:

\smallskip

\noindent $(IS')$ For every ultrafilter $u$ in $(B,\leq,\ll)$ such that $\exists c\in u$, $c\ll 1$, there exists a ultrafilter $v$ in $(B',\leq',\ll')$ such that $\exists c'\in v$, $c'\ll' 1$, and $\forall a\in u$ and $\forall b\in v$ there exists a $c_{ab}''\in B$ such that $c_{ab}''\ll 1$ and $\pl(b)\wedge c_{ab}''\not\ll(a\wedge c_{ab}'')^*$.
\end{defi}

\begin{theorem}\label{izdiskal}
The categories $\ISKAL$ and $\IISKAL$ are isomorphic; hence the categories $\SSKLC$ and $\IISKAL$ are dually equivalent.
\end{theorem}

\doc We will show that the restrictions $K_q:\ISKAL\lra \IISKAL$ and $\TE_q:\IISKAL\lra \ISKAL$ of the functors $K:\SKAL\lra \IKA$ and $\TE:\IKA\lra\SKAL$ defined in the proof of Theorem \ref{mvdeqvnk} are the desired isomorphism functors.

Let $\varphi\in\ISKAL((B,\rho,\BBBB),(B',\rho',\BBBB'))$. Since $\varphi$ is an $\IKA$-morphism (see the proof of Theorem \ref{mvdeqvnk}), we need only to check that $\varphi$ satisfies the axiom $(IS')$. Let $u$ be an ultrafilter in $(B,\leq,\ll)$ such that $\exists c\in u$, $c\ll 1$. Then, by \ref{eqv3}, $c\in\BBBB$ and hence $u$ is a bounded ultrafilter in $(B,\rho,\BBBB)$. Since $\varphi$ satisfies (IS), there exists a bounded ultrafilter $v$ in $(B',\rho',\BBBB')$, such that $\pl(v)\rho u$, i.e. $\exists c'\in v$, such that $c'\ll' 1$ and $\forall a\in u$ and $\forall b\in v$ there exists a $c_{ab}''\in B$ such that $c_{ab}''\ll 1$ and $\pl(b)\wedge c_{ab}''\not\ll(a\wedge c_{ab}'')^*$. Hence, $\varphi$ is an $\IISKAL$-morphism.

Let $\varphi\in\IISKAL((B,\leq,\ll),(B',\leq',\ll))$. We will show that the same function $\varphi:B\longrightarrow B'$ is a $\ISKAL$-morphism between $\TE_q(B,\leq,\ll)$ and $\TE_q(B',\leq',\ll)$. We have by \ref{mvdeqvnk} that $\varphi$ is an $\SKAL$-morphism. We need to check only that $\varphi$ satisfies $(IS)$. Let $u$ be a bounded ultrafilter in $(B,\rho,\BBBB)$. Then, by \ref{eqv3}, $\exists c\in u$, $c\ll 1$. Since $\varphi$ satisfies (IS'), there exists an ultrafilter $v$ in $(B',\leq',\ll')$ such that $\exists c'\in v$, $c'\ll' 1$ and $\forall a\in u$ and $\forall b\in v$ there exists an $c_{ab}''\in B$ such that $c_{ab}''\ll 1$ and $\pl(b)\wedge c_{ab}''\not\ll(a\wedge c_{ab}'')^*$. Hence $v$ is a bounded ultrafilter in $(B',\rho',\BBBB')$ and $\pl(v)\rho u$. Thus $\varphi$ is a $\ISKAL$-morphism. \sqs

In \cite{D-AMH2-10}, a category $\ISAL$ is introduced, namely, the objects of the category $\ISAL$ are all CLCA's (see \ref{locono}) and its morphisms are all injective complete Boolean homomorphisms $\varphi:(A,\rho,\BBBB)\longrightarrow(B,\eta,\BBBB')$ satisfying
axioms (L1)-(L3).

In \cite{D-AMH2-10}, a category $\SSAL$ is introduced, namely, the objects of the category $\SSAL$ are all CLC-algebras (see \ref{locono}) and its morphisms are all $\SAL$-morphisms $\p:(A,\rho,\BBBB)\lra
(B,\eta,\BBBB\ap)$ which satisfy condition (LS);

As it was proved in \cite{D-AMH2-10}, the categories $\SSLC$ and $\ISAL$ are dually equivalent. Also, in \cite{D-AMH2-10} it was proved that the categories $\ISLC$ and $\SSAL$ are dually equivalent.

\begin{defi}\label{soamorphin}
\rm
Let $\SIA$ be the category of all complete MVD-algebras and if $(B,\le,\ll)$ and
$(B\ap,\le\ap,\ll\ap)$ are two complete MVD-algebras then
$\p:(B,\le,\ll)\lra (B\ap,\le\ap,\ll\ap)$ will be an
$\SIA$-morphism iff $\p$ is an $\IA$-morphism (see \ref{catmvdn})
which satisfies the following axiom:

\smallskip

\noindent (CS) $\fa a,b\in B\ap$,  $a\ll\ap b$ implies $\pl(a)\ll
(\pl(b^*))^*$.
\end{defi}

\begin{defi}\label{convcatisn}
\rm We will denote by $\IIA$ the category of all complete
MVD-algebras and all injective complete Bool\-e\-an homomorphisms
between them satisfying  axioms (ES1), (S2) and (S3) (see
\ref{catmvdnk} and \ref{catmvdn}).
\end{defi}

\begin{theorem}\label{contheisn}
(i) The categories $\ISAL$ and $\IIA$ are isomorphic; hence the
categories $\SSLC$ and $\IIA$ are dually equivalent.

\noindent (ii) The categories $\SSAL$ and $\SIA$ are isomorphic;
hence the categories $\ISLC$ and $\SIA$ are dually equivalent.
\end{theorem}

\doc (i) It follows immediately from Theorem \ref{mvdeqvn}.

\medskip

\noindent (ii) Let
$\p\in\SSAL((B,\rho,\BBBB),(B\ap,\rho\ap,\BBBB\ap))$. We will
prove that the same function $\p:B\lra B\ap$ is an $\SIA$-morphism
between $K(B,\rho,\BBBB)$ and $K(B\ap,\rho\ap,\BBBB\ap)$ (see the
proof of \ref{mvdeqvnk} for $K$). In Theorem \ref{mvdeqvn} we have
seen that $\p$ is an $\IA$-morphism. Thus we need only to show
that $\p$ satisfies condition (CS) of \ref{soamorphin}. Put
$K(B,\rho,\BBBB)=(B,\le,\ll)$ and
$K(B\ap,\rho\ap,\BBBB\ap)=(B\ap,\le\ap,\ll\ap)$ (see \ref{mvdeqvn}
and \ref{eqv3} for the corresponding definitions). Let $a,b\in
B\ap$ and $a\ll\ap b$. Then, by (\ref{mvddefi3}),
$a(-\rho\ap)b^*$. Hence, by (LS), $\pl(a)(-\rho)\pl(b^*)$, i.e.
$\pl(a)\ll_\rho(\pl(b^*))^*$. Since, by (S2), $\pl(a)\ll 1$
(because $a\ll\ap 1$), we obtain, using twice (\ref{mvddefi3}),
that $\pl(a)\in\BBBB$ and $\pl(a)\ll (\pl(b^*))^*$. So, $\p$ is an
$\SIA$-morphism.

Let $\p\in\SIA((B,\le,\ll),(B\ap,\le\ap,\ll\ap))$. We  will show
that the same function $\p:B\lra B\ap$ is an $\SSAL$-morphism
between $\TE(B,\le,\ll)$ and $\TE(B\ap,\le\ap,\ll\ap)$ (see the
proof of \ref{mvdeqvnk} for $\TE$). For doing this it is enough
(by Theorem \ref{mvdeqvn}) to prove that $\p$ satisfies condition
(LS).
We will show that $\varphi$ satisfies condition (LS') which, as we know, is equivalent to condition (LS).

Put $\TE(B,\le,\ll)=(B,\rho,\BBBB)$ and
$\TE(B\ap,\le\ap,\ll\ap)=(B\ap,\rho\ap,\BBBB\ap)$.

Let $a,b\in B\ap$ and $a\ll_{\rho\ap} b$. Then, using
\ref{eqv3}, (S3), (\ref{L2rave}) and (CS), we obtain that $(a\llxa
b)\rightarrow (\fa c\in B$ such that $c\ll 1, \p(c)\we a\ll\ap
\p(c^*)\vee b)\rightarrow (\fa c\in B$ such that $c\ll 1,
\pl(\p(c)\we a)\ll (\pl(( \p(c^*)\vee b)^*))^*)\leftrightarrow
(\fa c\in B$ such that $c\ll 1,c\we\pl(a)\ll(\pl(\p(c)\we
b^*))^*)\leftrightarrow (\fa c\in B$ such that $c\ll
1,c\we\pl(a)\ll(c\we\pl(b^*))^*)\leftrightarrow (\fa c\in B$ such
that $c\ll 1,c\we\pl(a)\ll c^*\vee(\pl(b^*))^*)\leftrightarrow
(\pl(a)\llx(\pl(b^*))^*)$.

Therefore, $\p$ satisfies condition (LS). The rest follows from
Theorem \ref{mvdeqvn}.
\sqs

In \cite{D-AMH2-10}, a category $\IOPAL$ was introduced, namely, the objects of the category $\IOPAL$ are all CLCA's and its morphisms are all injective complete Boolean homomorphisms between them satisfying axioms (L1)-(L3) and (LO).
Also, in \cite{D-AMH2-10}, a category $\SOPAL$ was introduced, namely, the objects of the category $\SOPAL$ are all CLCA's and its morphisms are all surjective complete Boolean homomorphisms between them satisfying
axioms (L1)-(L3) and (LO).
As it was proved in \cite{D-AMH2-10}, the categories $\IOPAL$ and $\SOPAL$ are dually equivalent to the categories $\SOPLC$ and $\IOPLC$, respectively.

\begin{nota}\label{convcatis}
\rm We will denote by:

\noindent $\bullet$ $\IIPA$ the category of all complete
MVD-algebras and all injective complete Bool\-e\-an homomorphisms
between them satisfying  axioms (ES1), (S2), (S3) and (SO) (see
\ref{catmvdnk}, \ref{catmvdn} and \ref{catmvdo}).

\noindent $\bullet$ $\SIPA$ the category of all complete
MVD-algebras and all surjective complete Boolean homomorphisms
between them satisfying  axioms (ES1), (S2), (S3) and (SO) (see
\ref{catmvdnk}, \ref{catmvdn} and \ref{catmvdo}).
\end{nota}

\begin{theorem}\label{contheis}
(i) The categories $\IOPAL$ and $\IIPA$ are isomorphic; hence the
categories $\SOPLC$ and $\IIPA$ are dually equivalent.

(ii) The categories $\SOPAL$ and $\SIPA$ are isomorphic; hence the
categories $\IOPLC$ and $\SIPA$ are dually equivalent.
\end{theorem}

\doc It follows immediately from Theorems \ref{mvdeqv} and \ref{contheisn}. \sqs

In \cite{D-AMH2-10}, a category $\IOAL$ was introduced, namely, the objects of the category $\IOAL$ are all CLCA's and its morphisms are all complete Boolean homomorphisms between them satisfying axioms
(L1), (L2), (IS) and (LO);
Also, in \cite{D-AMH2-10}, a category $\SOAL$ was introduced, namely, the objects of the category $\SOAL$ are all CLCA's and its morphisms are all surjective complete Boolean homomorphisms between them satisfying axioms (L1), (L2) and (LO).

As it was proved in \cite{D-AMH2-10},
the categories $\IOLC$ and $\SOAL$ are dually equivalent; also, the categories $\SOLC$ and $\IOAL$ are dually equivalent.

\begin{nota}\label{i i s}
\rm We will denote by:

\noindent $\bullet$ $\IIOAL$ the category of all complete $MVD$-algebras and  all complete Boolean homomorphisms
between them satisfying the axioms $(S1)$, $(S2)$, $(IS')$ and $(SO)$ (see \ref{catmvdnk}, \ref{iiskal} and \ref{catmvdo}).

\noindent $\bullet$ $\ISOAL$ the category of all complete $MVD$-algebras and  all surjective complete Boolean homomorphisms between them satisfying the axioms $(S1)$, $(S2)$ и $(SO)$ (see \ref{catmvdnk} and \ref{catmvdo}).
\end{nota}

The next theorem follows immediately from Theorem \ref{mvdeqvnk}, Theorem \ref{mvdeqvo} and Theorem \ref{izdiskal}:

\begin{theorem}\label{iioal}
The categories $\IOAL$ and $\IIOAL$ are isomorphic; hence the categories $\SOLC$ and $\IIOAL$ are dually equivalent.
\end{theorem}

The next theorem follows immediately from Theorem \ref{mvdeqvnk} and Theorem \ref{mvdeqvo}:

\begin{theorem}\label{isoal}
The categories $\SOAL$ and $\ISOAL$ are isomorphic; hence the categories $\IOLC$ and $\ISOAL$ are dually equivalent.
\end{theorem}


\begin{defi}\label{conmvd3}
\rm An MVD-algebra $(B,\le, \ll)$ is called {\em connected}\/ if
it satisfies the following axiom:

\smallskip

\noindent (CONA) If $a\neq 0,1$ then there exists  $c\ll 1$ such
that $c\we a\nll a\vee c^*$.
\end{defi}

\begin{fact}\label{confact3}
Let $(L,\tau)$ be a locally compact Hausdorff space. Then the
standard MVD-algebra $(RC(L),\sbe,\ll_L)$ is connected iff the
space $(L,\tau)$ is connected.
\end{fact}

\doc Let's note that $\forall F,G\in RC(L)$, $F\ll_L G$ iff $F\in CR(L)$ and $F\ll_{\rho_L} G$. Hence $k(RC(L),\rho_L,CR(L))=(RC(L),\sbe,\ll_L)$. Let $(L,\tau)$ be connected. Then, from \ref{confact}, it follows that $(RC(L),\rho_L)$ is connected. Let $a\in RC(L)$, $a\not=0,1$. Then, from (CON), $a\rho_L a^*$. It follows from \ref{eqv3} that there exists $c\ll_L 1$ such that $c\we a\nll_L a\vee c^*$, i.e. $(RC(L),\sbe,\ll_L)$ is connected.

Let now $(RC(L),\sbe,\ll_L)$ be connected. Then for every $a\in RC(L)$ such that $a\not=0,1$, there exists an $c\ll_L 1$ such that $c\we a\nll_L a\vee c^*$. It follows from \ref{eqv3} that $a\rho_La^*$. Hence $(RC(L),\rho_L)$ is connected. Then, it follows from \ref{confact} that $(L,\tau)$ is connected. \sqs

\begin{nota}\label{IAC&IPAC}
\rm We will denote by:

\noindent $\bullet$ $\IAC$ the category of all connected complete
MVD-algebras and all complete Boolean homomorphisms between them
satisfying  axioms (ES1), (S2), (S3) (see \ref{catmvdn}).

\noindent $\bullet$ $\IPAC$ the category of all connected complete
MVD-algebras and all complete Boolean homomorphisms between them
satisfying  axioms (ES1), (S2), (S3) and (SO) (see \ref{catmvdnk},
\ref{catmvdn} and \ref{catmvdo}).
\end{nota}

The next theorem follows immediately  from \ref{confact},
\ref{confact3} and \ref{mvdeqvn}:

\begin{theorem}\label{conthen3}
The categories $\SALC$ and $\IAC$ are isomorphic; hence the
categories $\SLCC$ and $\IAC$ are dually equivalent.
\end{theorem}

The next theorem follows immediately from \ref{confact},
\ref{confact3}
and \ref{mvdeqv}:

\begin{theorem}\label{conthe3}
The categories $\OPALC$ and $\IPAC$ are isomorphic;
hence the categories $\OPLCC$ and $\IPAC$ are dually equivalent.
\end{theorem}

Analogously one can formulate and prove the connected versions of Theorems
\ref{mvdeqvnk} and \ref{mvdeqvo}.

In \cite{D-AMH1-10}, a category {\bf DHLC} was introduced, namely, the objects of the category {\bf DHLC} are all complete LC-algebras and its morphisms are all functions $\psi : (A,\rho,\mathbb{B})\longrightarrow(B,\eta,\mathbb{B}')$ between the objects of {\bf DHLC} satisfying the conditions

\smallskip

\noindent (DLC1) $\psi(0)=0$.

\noindent (DLC2) $\psi(a\wedge b)=\psi(a)\wedge\psi(b)$ for all $a,b\in A$.

\noindent (DLC3) If $a\in\mathbb{B}, b\in A$ and $a\ll_{\rho} b$, then $(\psi(a^*))^*\ll_{\eta}\psi(b)$.

\noindent (DLC4) For every $b\in\mathbb{B}'$ there exists $a\in\mathbb{B}$ such that $b\leq\psi(a)$.

\noindent (DLC5) $\psi(a)=\bv\{\psi(b)\st b\in\mathbb{B}, b\ll_{\rho}a\}$, for every $a\in A$.

Let the composition $``\odot$" of two morphisms $\psi_1:(A_1,\rho_1,\mathbb{B}_1)\longrightarrow (A_2,\rho_2,\mathbb{B}_2)$ and $\psi_2:(A_2,\rho_2,\mathbb{B}_2)\longrightarrow (A_3,\rho_3,\mathbb{B}_3)$ of {\bf DHLC} be defined by the formula

\begin{equation}\label{comp}
\psi_2\odot\psi_1(a)=\bigvee\{(\psi_2\circ\psi_1)(b)\st b\in\BBBB,b\ll_{\rho}a\},
\end{equation}

\noindent for every $a\in A$.

As it was proved in \cite{D-AMH1-10}, the category $\DHLC$ is dually equi\-valent to the category $\HLC$ of all locally compact Hausdorff spaces and all continuous mappings between them.

\begin{defi}\label{MVDHLC}
\rm
Let {\bf MVDHLC} be the category whose objects are all complete MVD-algebras and whose morphisms are all functions $\psi : (A,\leq,\ll)\longrightarrow(B,\leq',\ll')$ between the objects of {\bf MVDHLC} satisfying the conditions

\smallskip

\noindent (MVDLC1) $\psi(0)=0$.

\noindent (MVDLC2) $\psi(a\wedge b)=\psi(a)\wedge\psi(b)$ for all $a,b\in A$.

\noindent (MVDLC3) If $a,b\in A$ and $a\ll b$, then $\forall c\in B$ with $c\ll' 1$, $(\psi(a^*))^*\wedge c\ll'\psi(b)\vee c^*$.

\noindent (MVDLC4) For every $b\in B$ with $b\ll' 1$ there exists $a\in A$ with $a\ll 1$ such that $b\leq\psi(a)$.

\noindent (MVDLC5) $\psi(a)=\bv\{\psi(b)\st b\ll a\}$, for every $a\in A$.

Let the composition $``\circledcirc$" of two morphisms $\psi_1:(A_1,\leq_1,\ll_1)\longrightarrow (A_2,\leq_2,\ll_2)$ and $\psi_2:(A_2,\leq_2,\ll_2)\longrightarrow (A_3,\leq_3,\ll_3)$ of {\bf MVDHLC} be defined by the formula

\begin{equation}\label{comp1}
\psi_2\circledcirc\psi_1(a)=\bigvee\{\psi_2\circ\psi_1(b)\st b\ll_1 a\},\ \forall a\in A_1.
\end{equation}
\end{defi}

\begin{theorem}\label{lccont1}
The categories $\DHLC$ and\/ {\bf MVDHLC} are isomorphic; hence the categories {\bf MVDHLC} and {\bf HLC} are dually equivalent.
\end{theorem}

\doc Let us define two covariant functors $P:{\bf DHLC}\longrightarrow{\bf MVDHLC}$ and $Q:{\bf MVDHLC}\longrightarrow{\bf DHLC}$.

For every $(B,\rho,\BBBB)\in|{\bf DHLC}|$ we put $P(B,\rho,\BBBB)=k(B,\rho,\BBBB)$ (see \ref{eqv3} for $\k$). Then
Theorem \ref{eqv3} implies that $P$ is well-defined on the objects
of the category {\bf DHLC}.

Let $\psi\in{\bf DHLC}((B,\rho,\BBBB),(B',\rho',\BBBB'))$. We will prove that $\psi$ is a {\bf MVDHLC}-morphism between $P(B,\rho,\BBBB)=(B,\leq,\ll)$ and $P(B',\rho',\BBBB')=(B',\leq',\ll')$. It is obvious that $\psi$ satisfies axioms (MVDLC1) and (MVDLC2). Let $a\ll b$. Then $a\in\BBBB$ and $a\ll_{\rho}b$. It follows from (DLC3) that $(\psi(a^*))^*\ll_{\rho'}\psi(b)$. Then, from \ref{eqv3} it follows that $\forall c\ll' 1$, $(\psi(a^*))^*\wedge c\ll'\psi(b)\vee c^*$. Hence $\psi$ satisfies (MVDLC3).

Let $b\ll' 1$. From (DLC4) it follows that there exists an $a\in\BBBB$ such that $b\leq\psi(a)$. Hence $a\ll 1$ and $b\leq\psi(a)$, i.e. $\psi$ satisfies (MVDLC4).

Let $a\in B$. Then $\psi(a)=\bigvee\{\psi(b)\st b\in\BBBB, b\ll_{\rho}a\}=\bigvee\{\psi(b)\st b\ll a\}$. Hence $\psi$ satisfy (MVDLC5). Therefore $\psi\in{\bf MVDHLC}((B,\leq,\ll),(B',\leq',\ll'))$. So, we can define $P(\psi)=\psi$.

Let $\psi_i\in\DHLC((B_i,\rho_i,\BBBB_i),(B_{i+1},\rho_{i+1},\BBBB_{i+1}))$ and $P(\psi_i)=\varphi_i$, $i=1,2$. We have that $\forall a\in B_1$, $(\varphi_2\circledcirc\varphi_1)(a)=\bigvee\{(\varphi_2\circ\varphi_1)(b)\st b\ll_1 a\}=\bigvee\{(\psi_2\circ\psi_1)(b)\st b\in\BBBB_1, b\ll_{\rho_1}a\}=(\psi_2\odot\psi_1)(a)=(P(\psi_2\odot\psi_1))(a)$. Since, obviously, $P$ preserves the identities, we get that $P:{\bf DHLC}\longrightarrow{\bf MVDHLC}$ is a (covariant) functor.

Let $(B,\leq,\ll)\in|{\bf MVDHLC}|$. We put $Q(B,\leq,\ll)=\t(B,\leq,\ll)$ (see \ref{eqv3} for $\t$). Then
Theorem \ref{eqv3} implies that $Q$ is well-defined on the objects
of the category {\bf MVDHLC}.

Let $\psi\in{\bf MVDHLC}((B,\leq,\ll),(B',\leq',\ll'))$. We will prove that $\psi$ is a {\bf DHLC}-morphism between $Q(B,\leq,\ll)=(B,\rho,\BBBB)$ and $Q(B',\leq',\ll')=(B',\rho',\BBBB')$. It is obvious that $\psi$ satisfies axioms (DLC1) and (DLC2).

Let $a\in\BBBB$, $b\in B$ and $a\ll_{\rho}b$. Hence $a\ll b$. It follows from (MVDLC3) that $\forall c\ll' 1$, $(\psi(a^*))^*\wedge c\ll' \psi(b)\vee c^*$. By \ref{eqv3}, we get that $(\psi(a^*))^*\ll_{\rho'}\psi(b)$. Then $\psi$ satisfies (DLC3).

Let $b\in\BBBB'$. Then $b\ll' 1$. It follows from (MVDLC4) that $\exists a\ll 1$, such that $b\leq'\psi(a)$. Hence $a\in\BBBB$ and $b\leq'\psi(a)$. Then $\psi$ satisfies (DLC4).

Let $a\in B$. Then from (MVDLC5) we get that $\psi(a)=\bigvee\{\psi(b)\st b\ll a\}=\bigvee\{\psi(b)\st b\in\BBBB, b\ll_{\rho}a\}$. Therefore $\psi\in{\bf DHLC}((B,\rho,\BBBB),(B',\rho',\BBBB'))$. So, we can define $Q(\psi)=\psi$.

Let $\varphi_i\in {\bf MVDHLC}((B_i,\leq_i,\ll_i),(B_{i+1},\leq_{i+1},\ll_{i+1}))$, $Q(\varphi_i)=\psi_i$, $i=1,2$. We have that $\forall a\in B_1$, $(\psi_2\odot\psi_1)(a)=\bigvee\{(\psi_2\circ\psi_1)(b)\st b\in\BBBB_1, b\ll_{\rho_1}a\}=\bigvee\{(\varphi_2\circ\varphi_1)(b)\st b\ll_1 a\}=(\varphi_2\circledcirc\varphi_1)(a)=Q(\varphi_2\circledcirc\varphi_1)(a)$. Since, obviously, $Q$ preserves the identities, we get that $Q:{\bf MVDHLC}\longrightarrow{\bf DHLC}$ is a (covariant) functor.

From the definition of the functors $P$ and $Q$ and the equalities $\k\circ\t=id$, $\t\circ\k=id$ (see \ref{eqv3}), we conclude that $P\circ Q=Id_{\bf MVDHLC}$ and $Q\circ P=Id_{\bf DHLC}$. Hence, the categories {\bf DHLC} and {\bf MVDHLC} are isomorphic. \sqs

\baselineskip = 0.75\normalbaselineskip

\bibliographystyle{plain}

\end{document}